\documentclass[coll]{imsart}

\usepackage[square,numbers]{natbib}
\usepackage{amssymb}
\usepackage{amsthm,amsmath}
\usepackage{mathrsfs}
\usepackage{bm}
\usepackage{dsfont}
\usepackage{graphicx,color}

\allowdisplaybreaks

\startlocaldefs

\newcommand{\indic}{\mathds{1}}

\newcommand{\R}{\mathds{R}}
\newcommand{\N}{\mathds{N}}

\newcommand{\Bcr}{\mathscr{B}}

\newcommand{\Ecr}{\mathscr{E}}

\newcommand{\Lcr}{\mathscr{L}}

\newcommand{\ddr}{\mathrm{d}}

\newcommand{\idr}{\mathrm{i}}
\newcommand{\Log}{\mathrm{Log}\,}

\newcommand{\ep}{\varepsilon}

\endlocaldefs

\begin{document}

\begin{frontmatter}

\title{The origins of de Finetti's critique of countable additivity}
\runtitle{de Finetti and countable additivity}


\author{\fnms{Eugenio} \snm{Regazzini}\thanksref{t1}\ead[label=e1]{eugenio.regazzini@unipv.it}}
\address{
Dipartimento di Matematica ``F. Casorati'' 
Universit\`a di Pavia (Italy)
\printead{e1}}

\runauthor{E. Regazzini}

\thankstext{t1}{Also affiliated to CNR-IMATI, Milano (Italy). This work is partially supported by the Italian Ministry of University and Research, grant n.~2008MK3AFZ. The author also wishes to thank Antonio Lijoi and an anonymous Referee for their  comments that led to an improved exposition.} 

\affiliation{University of Pavia}
\begin{abstract}
Bruno de Finetti was one of the most convinced advocates of finitely additive probabilities. The present work describes the intellectual process that led him to support that stance and provides a detailed account both of the first paper by de~Finetti on the subject and of the ensuing correspondence with Maurice Fr\'echet. Moreover, the analysis is supplemented by a useful picture of  de~Finetti's interactions with the international scientific community at that time, when he elaborated his subjectivistic conception of probability.
\end{abstract}

\begin{keyword}[class=AMS]
60-03, 60A05
\end{keyword}

\begin{keyword}
Finitely additive probabilities; Processes with independent and
stationary increments; Subjective probabilities
\end{keyword}

\end{frontmatter}

\section{Introduction}

Finitely additive probabilities  are indissolubly linked to the name of Bruno de~Finetti (1906-1985). Indeed, he has been one of the most convinced advocates of finite additivity, since he has started working on the mathematical formulation that he proposed, in 1930, for his subjectivistic conception of probability. Most of the recent contributions to this topic in the literature rely on (English translations of) late works by de~Finetti, instead of considering his early papers containing a wealth of fresh and original ideas. A typical feature of de~Finetti's late works is that they generally aim at providing critical syntheses of his original way of thinking about crucial problems concerning: Probability, Induction, Statistics, Insurance, Economics, Politics, to say nothing of the philosophical debate at large. Hence, in these writings he makes limited use of mathematical formalism, omits precise references of an historical nature, whilst he often jumps at the chance of both controversial amusing hints and sharp provocative cues. Nevertheless, late works share with the early ones the feature of being faultless from the point of view of logical and conceptual accuracy. It is not, then, surprising that, as far as finite additivity is concerned, essays on de~Finetti's work often reduce the topic either to an intellectual activity or, at best, to an issue of a mere mathematical taste. In the latter case, special attention is given to connections and consequences of a formal nature. So, they generally fail to shed some light on the intellectual efforts that led him to share that seemingly singular position. In fact, the authentic motivations supporting de~Finetti's stance in probability can be found in his early works. These, among other things, are tersely written and contain accurate formulations of theorems thus allowing one a more sound understanding of his innovations.

The main purpose of the present work is to discuss the first paper that, at the best of this author's knowledge, de~Finetti devoted to the analysis of the effects of considering as admissible finitely additive laws. See \cite{df(30a)}. The title of the paper, ``\textit{Sui passaggi al limite nel calcolo delle probabilit\`a}'', could be translated into ``\textit{On the limit processes in the calculus of probability}'' and evokes the continuity property of countably additive laws. Its content is carefully described and critically annotated in Section~4. This follows a discussion, in Section2, on the development of de~Finetti's ideas and achievements in probability theory between 1927, the year of his graduation, and 1930, the year of the publication of ``\textit{Sui passaggi~...}''. For a more comprehensive illustration of de~Finetti's work, see \cite{cr(96)}. Section~3 provides a concise description of the mathematical theory of probability deduced from de~Finetti's \textit{coherence principle}. See \cite{df(31b)}.

The present author is not aware of how ``\textit{Sui passaggi~...}'' was received by the Italian scientific community. However, the fact that it was not published in the journal \textit{Rendiconti della Reale Accademia dei Lincei} might conceal a cold reaction by the most influential Italian probabilist of the time, namely Francesco Paolo Cantelli (1875-1966). It is well-known that Cantelli was not enthusiastic about non countably additive probabilities. See \cite{cant(17)} and \cite{cant(35a), cant(35b), cant(36)}. A public reaction on the spur of the moment came from the famous French mathematician Maurice Fr\'echet (1878-1973) who initiated an interesting correspondence with de~Finetti gathered into four published notes. See \cite{df(30b), df(30c)}, \cite{f(30a), f(30b)}. The reading of this correspondence provides a vivid insight into the stances of the two ``competitors'' that reflect two different ways of thinking still of great interest. Then, this correspondence is reported and annotated in the present paper as well. See Section~5.

I am delighted I am given the chance to write this work in honour of Joe Eaton who devoted a distinguished part of his scientific research to the foundations of probability and statistics.

\section{De Finetti and other probabilists at the end of the twenties}

De Finetti started coping with the fundamental task of formulating a satisfactory and general theory of probability right after his graduation. Almost in the same period he initiated his studies on the sequences of exchangeable events and on functions with independent, stationary increments (f.i.s.i., for short). But, whilst his investigations on these specific topics proceeded rapidly~-- cf., e.g., \cite{r(07)} and \cite{br(08)}~-- the research on the foundations of probability appeared somewhat bristling with difficulties. The major challenge arose from the mathematical formulation of the subjectivistic conception of probability as expounded, from a philosophical viewpoint, in \cite{df(31a)}, an essay that had already been drawn up at the beginning of 1930. That time lag did not restrain him from reaching considerable achievements, but led him to a retrospective critical reflection on some of his own results, that he had obtained within the realm of countably additive probabilities. This happened during the second half of 1929, when he reached the conclusion that such a condition is not necessary in order that a function on a class of events can be viewed as a probability. Indeed the presentations, at the \textit{R. Accademia dei Lincei}, of \cite{df(29b)} and \cite{df(29c)} are the 7th of October and the 1st of December, respectively. The precise deduction of the characterizing properties of a probability law, from the subjectivistic standpoint, were announced in both the replies to Fr\'echet. See \cite{df(30b), df(30c)} and \cite{df(31b)} for the final version.

It is worth providing some further insight into this story. By assuming the continuity  of the law of a f.i.s.i. based on Gaussian finite--dimensional distributions, in \cite{df(29a)} de~Finetti had stated that almost every trajectory turns out to be nowhere differentiable: a result today ascribed to Paley, Wiener and Zygmund. See \cite{pwz(33)}. Subsequently, he had shown that the probability distribution (p.d.) of $X_t$ is continuous, for each $t>0$, whenever $(X_t)_{t\ge 0}$ is a continuous f.i.s.i.. See \cite{df(29b)}. A first clear statement of his concerns about the suitability of an analysis confined to considering random functions governed by countably additive laws can be found in \cite{df(29c)} where he deals with the problem of determining the p.d. of the integral, on $[0,t]$, of a f.i.s.i.. See next Section~4. It is in this last essay that he announces the preparation of \textit{Sui Passaggi}, with the purpose of investigating the issue from a more general viewpoint. It seems fair to say that \textit{Sui Passaggi} marks a turning point in de~Finetti's mathematical treatment of probability. After hinting at the possible effects of such an afterthought on his own work, we linger on describing the boundary conditions within which it matured.

As to the Italian scientific community, it has been already mentioned that Cantelli claimed to be unconditionally in favour of countable additivity. As an example, in his celebrated paper on the strong law of large numbers~-- see \cite{cant(17)}~-- he had declared ``Such an assumption plainly cannot raise objections from a theoretical viewpoint, and corresponds to the feeling that probability, viewed from an empirical angle, arouses in us''. He made no mention of the fundamental problem of the existence of a countably additive extension of sequences of assigned finite--dimensional laws of $(X_1,\ldots,X_n)$, for $n=1,2,\ldots$, to obtain a probability law for $(X_n)_{n\ge 1}$. See also \cite{s(92)} and \cite{r(05)}. We do not know whether he ignored or not the existence of such a problem that was completely solved later by Kolmogorov. See \cite{k(33)}. In fact, it was only in 1932 that Cantelli proposed a measure--theoretic approach with an explicit view to proving existence of extensions like the previous one. See, in particular, \cite{o(39)}, which includes one of the main achievements within the so--called \textit{Cantelli abstract theory of probability}. As far as the issue of interpreting the meaning of probability is concerned, Cantelli swung between an empirical interpretation and a bent for finding conceptual connections between his abstract theory and the classical Laplace definition. It is for sure that he had no sympathy for subjectivistic interpretations.

As regards his interactions with the broader international scientific community, de~Finetti has been keeping up correspondences with some of the foreign foremost scholars such as Fr\'echet, L\'evy, Khinchin and Kolmogorov since the end of the 20's. It is well--known, for example, that Alexandr Y. Khinchin (1894--1959) took interest in the study of sequences of exchangeable events during and after the \textit{8th International Congress of Mathematicians} held in Bologna, 3--10 September, 1928. See \cite{kh(32a), kh(32b)}. Andrei N. Kolmogorov (1903-1987), following the theory devised by de~Finetti in '29, obtained the renowned representation of the characteristic function of $X_1$, when $(X_t)_{t\ge 0}$ is a f.i.s.i. and $X_1$ has finite variance. The starting point for this research is in \cite{df(30d)}, whereas the Kolmogorov contribution is contained in \cite{k(32a), k(32b)}. Also the generalization due to Khinchin of the Kolmogorov statement can be viewed, unlike that given by L\'evy, as a by--product of the de~Finetti approach to functions with independent increments. See \cite{kh(37)} and \cite{l(34)}. In the period 1929-1930, de~Finetti focused on a paper by Kolmogorov concerning the representation of associative means. See \cite{k(30)} and \cite{df(31c)}. Both Fr\'echet and Evgeny E. Slutsky (1880--1948) were mentioned by de~Finetti for contributions to stochastic convergence: the former for the study of convergence of random elements in abstract spaces (see \cite{df(30b)}), the latter for the use of the term \textit{stochastic} that de~Finetti adopted to designate convergence in probability (see \cite{df(30a)}). As to Paul L\'evy (1886--1971), it is well--known that he conceived and developed his fundamental contribution to f.i.s.i.'s independently of both de~Finetti and Kolmogorov. See \cite{l(70)}. On de~Finetti's side, L\'evy's \textit{Calcul des Probabilit\'es} (1925) and Castelnuovo's \textit{Calcolo delle Probabilit\`a e Applicazioni} (1919) were the sole existing reference books at the end of the twenties, gathering the essentials of the theory and calculus of probability in a systematic way. Like Cantelli, L\'evy was inclined to restrict probability laws to countably additive functions on events, and tried to justify his position in the final part of the book.

As to the meaning and the interpretation of the concept of probability, the span of time we are dealing with overlaps with the success of the so--called empirical conceptions, according to which probability is related to frequency. This position was defended, in those years, by great scientists like Guido Castelnuovo (1865--1952) and Richard von Mises (1883--1953). See, for example, \cite{vm(28)}. De Finetti  had found the empirical arguments incomplete and inadequate since the very beginning of his own approaching the probabilistic studies. Unsatisfied, he followed a different path that led him to the formulation of a radically subjectivistic theory, with the consequence that countable additivity is not necessary in order that a set function may be considered a probability. A brief description of de~Finetti's theory will be sketched in the next section. We conclude the present one with a digression about the spread of basic elementary concepts in the literature of the day, by means of an example. The research into the subject has been suggested by the reading of \textit{Sui Passaggi} and regards, in particular, definition and properties of convergence in probability.

In Section~8 of that paper, de~Finetti reminds the reader of a definition given in \cite{df(29c)}~-- the definition of \textit{convergence in probability} of sequences of random numbers~-- and proves that convergence in probability entails \textit{convergence in distribution}. It will be explained, in Section~4, why he had been led to deal with this arrangement of topics. Here we want to comment on the fact that in neither of the two papers does he mention any reference and he proposes to designate the concept by the locution \textit{convergenza stocasticamente uniforme} (stochastically uniform convergence). Combining this circumstance with the fact that de~Finetti was used to carefully mention references, including those having little bearing on the development of his own research, one is led to conclude that the concept of convergence in probability was not yet a part  of the probabilistic literature at the time. This may seem amazing. It is worth recalling that the topic of convergence of sequences of random numbers had generated heated controversies, not yet dulled at the end of the twenties. A significant part of the contention can be traced back to Cantelli's determination in claiming his priority, over Emile Borel (1871--1956), for the formulation of the strong law of large numbers. See, for example, \cite{s(92)} and references therein. Stimulated by the draft of the present paper, we have thoroughly examined the literature of that time and we have got to the conclusion that de~Finetti was probably the first to consider sequences having a random number as a limit. Moreover, one should acknowledge his priority in proving that convergence in probability is stronger than convergence in distribution. The statement on page 25 of \cite{s(92)}, according to which ``Cantelli (1916) had anticipated Slutsky (1925) in introducing a random variable (rather than just a constant) as a limit ... of convergence in probability'', does not seem correct. Indeed, on the one hand, Cantelli in \cite{cant(16)} deals with constant limits only~--see also \cite{cant(35a)}~-- and, on the other hand, Slutsky in \cite{sl(25)} considers only sequences of the type $(X_n-m_n)_{n\ge 1}$ 
converging to zero for some sequence $(m_n)_{n\ge 1}$ of real numbers (\textit{stochastische Asymptote}). A more delicate analysis must be devoted to the Fr\'echet work on convergence of random sequences. It is Fr\'echet himself who gave notice, in his comments on \textit{Sui Passaggi} that were presented on the 3rd of July 1930 at the \textit{Reale Istituto Lombardo}, of a work of his having some points of contact with de~Finetti's paper. He did not mention the title of his article, which was referred to as a summary of his recent courses, and announced it was going to appear in the (Italian) statistics journal \textit{Metron}. Indeed, it was published in the last issue of that year. See \cite{f(30a)} and \cite{f(30c)}. In the meantime, de~Finetti has come by the Fr\'echet paper: this is witnessed by the fact he touches upon it in his reply to the Fr\'echet first group of comments on \textit{Sui Passaggi}. De~Finetti confines himself to saying that, apart from what is consequence of the assumption of countable additivity, he has particularly appreciated the study of the convergence of sequences of random elements taking values in abstract metric spaces. See the last part of \cite{f(30c)}. The first sections contain accurate descriptions of the concepts of convergence in probability and of almost sure convergence for sequences of random numbers. As to the former, the definition is obviously the very same as that previously proposed, under different name, by de~Finetti. Moreover, also the Fr\'echet work contains a proof of the fact that convergence in probability entails convergence in distribution. The funny thing is that he omits giving precise references to the points of contact with \textit{Sui Passaggi}. In fact, he just mentions works by Slutsky, Dell'Agnola and de~Finetti, without specifying either titles or other bibliographical data. He justifies these omissions by writing that ``... nous n'avons pu consulter assez librement et compl\`etement les publications de M. Slutsky ... Par contre les m\'emoires de M.M. Dell'Agnola et de~Finetti ne nous paraissent pas consacr\'es aux m\^emes questions que celles qui ont \'et\'e trait\'ees ici.'' Such a comment 
is not appropriate either for the Dell'Agnola paper \cite{da(29)} or for \textit{Sui Passaggi}. But, while the former will be mentioned in Chapter~5 of \cite{f(37)}, de~Finetti's contribution will be ignored even on that occasion, in spite of the anything but negligible overlapping. So, de~Finetti was probably the first to introduce a general definition of  convergence in probability and to study its connections with convergence in law, but Fr\'echet was closely following him on this path. Situations of this type are not infrequent in de~Finetti's scientific production. Other noteworthy examples are represented by: (i) the concept of inifinitely divisble law, introduced to encapsulate the characterizing feature of the law of $X_1$ in a stochastic f.i.s.i. $(X_t)_{t\ge 0}$ (see \cite{df(29a), df(30d)}); (ii) the completion of the continuity theorem for characteristic functions, to provide the first proof of the renowned representation theorem for exchangeable events (see \cite{df(30f)}); (iii) the extension of the  Kolmogorov theorem for associative means (see \cite{df(31c)}); (iv) the anticipation of the so--called Kendall's $\tau$ coefficient, in a general study on correlation and concordance (see \cite{df(37)}); (v) a pioneering use, almost contemporaneously with L\'evy, of the notions of martingale, stopping time and optional sampling, to reformulate the Lundberg--Cram\'er theory of risk, without coining any new locution to designate them (see \cite{df(39)}). The recurrence of these circumstances is the result of a combination of his special mathematical inventiveness with the fact that, as he wrote himself, ``he was interested in mathematics meant as a tool for applications ... and for the investigation of conceptual and critical issues ..., rather than as formalism or as an abstract subject, axiomatized and withdrawn in itself.'' See Page XVIII of \cite{df(81)}. This attitude prevented him from isolating results that were merely proved with a view to the solution of a more general problem. It also made it difficult to acknowledge the paternity of a number of scientific innovations.

\section{De Finetti's coherence principle}

Responding to Fr\'echet's second series of remarks, de~Finetti claims that the most significant point of the question on additivity lies in the need of proofs for the properties (of probability) that one wishes to affirm. See \cite{df(30c)}. For those who, like us, are accostumed to affirm those properties by means of axioms, de~Finetti's recommendation might appear as the fruit of an outdated way of thinking. However, it is plain to see that it represents the most reasonable way to settle the dispute constructively. Then, he asks every author to provide a proof of one's own, consistent with a well-specified conceptual starting point. For a complete understanding of the coming sections, it is useful to recall the solution to which de~Finetti himself always made reference, starting from the end of 1929. See \cite{df(30e)} and \cite{df(31b)}.

He maintains probability assessments definitely have a subjectivistic source, but one can conceive the following ideal experiment to check on the closeness of a real--valued function $P$, defined on a class $\Ecr$ of events, to one's authentic opinions on the uncertainty of the elements of $\Ecr$. One ought to think of $\{P(E):\: E\in\Ecr\}$ as a system of unit prices to have bets on the events included in $\Ecr$. More precisely, if one is willing to accept bets of any amount (either positive or negative)  at the unit prices fixed by the above system, then $P$ represents the desired quantification of one's judgments. De~Finetti assumes an ounce of rationality, in the sense that $P$ is required to ensure that there is no choice of $c_1,\ldots,c_n$ in $\R$ and of $E_1,\ldots,E_n$ in $\Ecr$ such that $\sup \sum_{k=1}^n c_k\{P(E_k)-\indic_{E_k}\}<0$, where $\sup$ is taken with respect to all the elements (elementary events) of the partition of the sure event $\Omega$ generated by $\{E_1,\ldots,E_n\}$, and $\indic_E$ denotes the indicator function of the event $E$. Indeed, with this notation, $\sum_{k=1}^n c_k\{P(E_k)-\indic_{E_k}\}$ represents the gain from a combination of bets of amount $c_1,\ldots,c_n$ on $E_1,\ldots,E_n$, respectively. These remarks led de~Finetti to restrict the class of admissible probability laws, on $\Ecr$, to those which obey the \textit{coherence principle}, i.e.:
\smallskip

\noindent For all \textit{finite} subsets $\{E_i:\: i\in I\}$ of $\Ecr$ and $\{c_i:\: i\in I\}$ of $\R$, $P$ satisfies
  \begin{equation}
    \label{eq:coehrence}
    \sup \sum_{i\in I} c_i\{P(E_i)-\indic_{E_i}\}\ge 0
  \end{equation}
with $\sup$ taken with respect to all elementary events relative to $\{E_i:\: i\in I\}$.
\smallskip

Any $P$ satisfying this property is said to be a \textit{probability} on $\Ecr$.

Existence of at least one probability law on $\Ecr$ is proved in \cite{df(30e)}. It is of great interest the fact that coherence is well-defined, independently of the structure of $\Ecr$. Moreover, as proved in \cite{df(30e)}, any probability on $\Ecr$ admits a coherent extension to any larger class of events.

As to the properties that one wishes to affirm, with a view to the calculus of probabilities, it is easy to \textit{prove} that coherence entails:
\begin{itemize}
\item[($\pi_1$)] If $\Omega\in\Ecr$, then $P(\Omega)=1$
\item[($\pi_2$)] $P(E)\ge 0$ for every $E$ in $\Ecr$
\item[($\pi_1$)] If $E_1$, $E_2$ and $E_1\cup E_2$ are in $\Ecr$, with $E_1\cap E_2=\varnothing$, then
\[
P(E_1\cup E_2)=P(E_1)+P(E_2).
\]
\end{itemize}
Moreover: If $\Ecr$ is an algebra, then $(\pi_1)$, $(\pi_2)$ and $(\pi_3)$ are also sufficient in order that $P:\Ecr\to\R$ can be considered a probability. See \cite{df(31b)}.

That is de Finetti's \textit{proof} that countable additivity is not necessary for a function to be a probability on an algebra of events.

It would be misleading to consider countable families of bets with the purpose of extending additivity to countable families of pairwise disjoint events. Indeed, if one proceeds in this direction, it would be necessary to introduce extra--conditions, completely extraneous to the common interpretation of the term probability, in order to give a precise and unambiguous meaning to the gain.

To conclude, we mention a couple of facts that are of importance for understanding a few of the coming examples. Let $(P_n)_{n\ge 1}$ be a sequence of probabilities on $\Ecr$ and let $\Lcr:=\{E\in\Ecr:\: \lim_{n\to\infty}P_n(E) \mbox{ exists}\}\ne \varnothing$. Then $Q(E):=\lim_{n\to\infty}P_n(E)$ is a probability on $\Lcr$. It should be noted that $Q$ is not necessarily countably additive even if each $P_n$ is countably additive. The second fact to be considered herein is concerned with the general theory of stochastic processes. De~Finetti has frequently referred to it, even if tacitly, in his already mentioned papers on random functions. Only at the end of the last century it was noticed and encapsulated into a theorem by Lester E. Dubins (1920--2010). See \cite{dub(99)}. Define two (real--valued) random functions on $[0,+\infty)$ to be \textit{cousins} if the family $J$ of finite--dimensional p.d.'s of one of the functions is the same as the $J$ of the other random function. Dubins proves that each random function, in particular the Brownian motion, has a cousin almost all of whose paths are polynomials, another cousin almost all of whose paths are step functions (on each bounded time--interval, they only have a finite number of values, each assumed on a finite union of intervals) that are continuous on the right (on the left), and a fourth cousin almost all of whose paths are continuous, piecewise--linear functions. Hence, in the next sections there will be no contradiction when, referring to de~Finetti, continuous random functions will be considered with independent and stationary increments, different from the Brownian motion.

\section{Presentation and critical comment of \textit{Sui Passaggi}}

The paper we are now going to analyze consists of nine sections. The first three are devoted to the explanation of some elementary facts concerning p.d. functions: these provide simple illustrations of how certain conclusions, valid for countably additive probabilities, do not generally hold any more for finitely additive probabilities. Sections 4 to 8 deal with convergence of sequences of random variables and include interesting , and somewhat amazing, remarks about the Cantelli strong law of large numbers. Finally, in Section~9, de~Finetti goes back over the problem that had led him to the reflection developed in the previous sections: Is the p.d. of the integral of a random function equal to the limit of the laws of the integral sums?

\subsection{Discussion of Sections 1--3} De Finetti starts with the p.d. function of a random number $X$, defined as
\[
F(x)=P\{X < x\}\qquad\qquad\quad (x\in\R)
\]
and lingers on the correct interpretation of the discontinuity points of $F$, within a finitely additive framework. In order to avoid unnecessary complications, think of $P$ as a probability on $2^\Omega$, $\Omega$ being a set and $X$ a real--valued function defined on $\Omega$. See the previous section for the definition and the existence of a probability on $2^\Omega$, and for de~Finetti's attitude with respect to the general theory of probability at the end of 1929, when he was writing \textit{Sui Passaggi}. In particular, he was perfectly aware of the existence of probabilities $P$ such that $P(D)=1$ and $P(\{x\})=0$ for some countable $D\subset\Omega$ and for every $x$ in $D$. Since this will occur quite often in discussing de~Finetti's paper, we now indicate a way to obtain probabilities of that type.
\medskip

\textbf{Example 4.1.} Let $D$ be a countably infinite subset of $\Omega$, say $D=\{x_1,x_2,\ldots\}$, and let $P_n$ be a probability on $2^\Omega$ defined by
\[
P_n(A):=\frac{\sharp (A\cap D_n)}{n} \qquad\qquad (A\subset\Omega)
\]
with $D_n:=\{x_1,\ldots,x_n\}$. Write $\Lcr:=\{A\subset\Omega:\: \lim_{n\to\infty}
\sharp(A\cap D_n)/n \mbox{ exists}\}$ and set $P(A):=\lim_{n\to\infty}P_n(A)$ for every $A$ in $\Lcr$. As recalled in the previous section, $P$ is a probability on $\Lcr$ and, then, it admits a coherent extension $\gamma$ on $2^\Omega$ such that $\gamma(\{x_k\})=0$ and $\gamma(D)=1$. \hfill\qed
\medskip

We are in a position to discuss the main issue of the present section, i.e. the correct interpretation of discontinuity for a p.d. function, within the frame established in Section~3. Under the ordinary condition of countable additivity, for any $d$ of such type one would get $F(d-0)=F(d)=P\{X<d\}<F(d+0)=P\{X\le d\}$ and, consequently, $F(d+0)-F(d-0)$ would represent the probability \textit{concentrated} in $\{d\}$. This statement is not necessarily true if $P$ is just finitely additive: in this case, one can only say that the following chain of inequalities holds true:
\begin{equation}
  \label{eq:inequal}
  F(d-0)\le P\{X<d\} \le P\{X\le d\}\le F(d+0)
\end{equation}
along with $F(d-0)<F(d+0)$ if $d$ is a discontinuity point.

What de Finetti points out as a frequent mistake~-- in which, as recalled in Section~2, he himself had been trapped~-- the fact that many authors resorted to passages to the limit without preventive proof of the necessity of any condition justifying the limit processes along monotone sequences of events. Assuming either of the equalities $F(d-0)=P\{X<d\}$, $F(d+0)=P\{X\le d\}$ might be an istance of that mistake, as displayed in the following example drawn from \textit{Sui Passaggi}.
\medskip

\textbf{Example 4.2.} Let $\Omega=\R$ and $(x_n)_{n\ge 1}$ be a sequence with $x_{n+1}<x_n$, for every $n$, and $x_n\downarrow 0$. Define $D$ and $\gamma$ as in the previous Example~4.1, and the random number $X$ by $X(\omega)=\omega$ for every $\omega$ in $\Omega=\R$.  Then, the p.d. function $F$ of $X$ satisfies $F(x)=1-P\{X\ge x\}=1$ for every $x>0$, and $F(x)=P\{X\le x\}=0$ for every $x\le 0$. So, the jump $(=1)$ of $F$ at $d=0$ does not represent a mass concentrated in $\{0\}$. \hfill \qed
\medskip

This example shows that the inequality to the right of \eqref{eq:inequal} cannot be replaced by equality, excepting special cases. To see that an analogous remark can be made for the inequality to the left, it is enough to consider the case of a sequence $(y_n)_{n\ge 1}$ with $y_n<y_{n+1}$ for every $n$ and $y_n\uparrow 0$. Now, replacing $P$ with a probability $Q$ for which $Q(\{y_k\})=0$ for every $k$ and $Q(\{y_1,y_2,\ldots,\})=1$, the consequent p.d. function $G$ of $X$ satisfies $G(x)=0$ for $x<0$ and $G(0)=1=P\{X<0\}$. It should also be noted that $H:=pF+(1-p)G$ is a p.d. function for each $p$ in $[0,1]$, and for every $p$ in $(0,1)$ inequalities under discussion hold on both sides of \eqref{eq:inequal} when $F$ is replaced by $H$ and $P$ by $pP+(1-p)Q$. It is also straightforward to find variants exhibiting concentrated masses, which satisfy
\[
P\{X=d\}< F(d+0)-F(d-0).
\]
As to the behaviour of a p.d. function $F$ at $-\infty$ ($+\infty$, respectively), what can be said, in general terms, is that $\lim_{x\to -\infty}F(x)\ge 0$ ($\lim_{x\to+\infty} F(x)\le 1$, respectively), strict inequalities being possible,

\subsection{Discussion of Sections 4--9}

In Sections 4--6, de~Finetti discusses an important case in which, on the basis of the tacit assumption of continuity for probabilities, authors of the time were led to endow strong laws of large numbers with a meaning more general than one's due. De~Finetti focuses on Cantelli's proof of the convergence of the frequency of success relative to $n$ trials, in a sequence of Bernoulli trials, as $n$ goes to infinity. To tackle the problem in the usual terms, define $\Omega$ to be the set of all sequences $d:=(d_n)_{n\ge 1}$, each $d_n$ being $0$ or $1$. Then, for each $n$ define the $n$--th projection $p_n(d):=d_n$ $(d\in\Omega)$, and set
\[
f_n(d)=\frac{1}{n}\sum_{j=1}^n p_j(d).
\]
Then, $f_n$ represents the frequency of $1$ in the first $n$ trials. Let $P$ be any probability on $2^\Omega$ such that
\[
P(\{d\in\Omega:\: p_1(d)=e_1,\ldots,p_n(d)=e_n\})=p^{\sum_{i=1}^n e_i}
(1-p)^{n-\sum_{i=1}^n e_i}
\]
where $(e_1,\ldots,e_n)$ is any point in $\{0,1\}^n$, with $n=1,2,\ldots$, $p$ some fixed point in $[0,1]$, and provided that $0^0= 1$. Under these conditions, which imply that $(p_n)_{n\ge 1}$ is a Bernoulli sequence, Cantelli had proved that, for every $\ep,\delta>0$, there is $n_0=n_0(\ep,\delta)$ such that
\begin{equation}
  \label{eq:conv_cantelli}
  \inf_{m\ge 1} P\left(\left\{\max_{n\le k\le n+m}|f_k-p|\le \ep\right\}\right)\ge 1-\delta
\end{equation}
holds for any $n\ge n_0$.

It must be said that assessments of $P$ for which
\[
\inf_{m\ge 1} P\left(\left\{\max_{n\le k\le n+m}|f_k-p|\le \ep\right\}\right)= P\left(\left\{\max_{ k\ge n}|f_k-p|\le \ep\right\}\right),
\]
obtained by interchanging $\lim_{m\to\infty}$ with $P$, are consistent with the Bernoulli assumptions even in the frame of de~Finetti's theory, but they are not the sole. In particular, there are assessments for which one cannot say that $(f_n)_{n\ge 1}$ converges to $p$ with probability one, even though \eqref{eq:conv_cantelli} is obviously valid. In Section~5 de~Finetti illustrates the situation by means of the following interesting example.
\medskip

\textbf{Example 4.3.} In the space $\Omega$ of all sequences $d$, each $d_n$ being $0$ or $1$, define $S_n$ to be the set of all sequences $(e_1,\ldots,e_{n-1},1,0,0,\ldots)$ obtained, for each $n\ge 2$, as $(e_1,\ldots,e_{n-1})$ varies in $\{0,1\}^{n-1}$. Moreover, set $S_1:=\{1,0,0,\ldots\}$. So $\Omega_1:=\cup_{n\ge 1}S_n$ is the set of all sequences in $\Omega$ with a place occupied by $1$ and followed by an entire sequence of $0$'s. Now, consider the sequence $(Q_n)_{n\ge 1}$ of probabilities on $2^\Omega$ defined by $Q_1(S_1)=1$ and, for any $n\ge 2$
\[
Q_n(\{e_1,\ldots,e_{n-1},1,0,0,\ldots\})=p^{\sum_{i=1}^{n-1}e_i}
(1-p)^{n-1-\sum_{i=1}^{n-1}e_i}
\]
for any $(e_1,\ldots,e_{n-1})$ in $\{0,1\}^{n-1}$. Clearly $Q_n(S_n)=1$ for every $n\ge 1$. Finally, consider the probability $\gamma$ of Example~4.1, with $\Omega=\N:=\{1,2,\ldots\}$ and, for any subset $A$ of $\Omega$, set
\[
Q(A):=\int_\N Q_n(A)\,\gamma(\ddr n)
\]
the integral being meant in the sense of Riemann-Stieltjes as explained, for example, in \cite{br(83)}. It is easy to verify that $Q_n(\{p_1=e_1,\ldots,p_N=e_N\})=p^{\sum_{i=1}^N e_i}(1-p)^{N-\sum_{i=1}^N e_i}$ holds true for every $n\ge N+1$ and, since $\gamma(\{N+1,N+2,\ldots\})=1$ for every $N$, one gets
\begin{align*}
Q(\{p_1=e_1,\ldots,p_N=e_N\})
&=\int_{\{n\ge N+1\}}Q_n(\{p_1=e_1,\ldots,p_N=e_N\})\,
\gamma(\ddr n)\\[7pt]
&=p^{\sum_{i=1}^N e_i}(1-p)^{N-\sum_{i=1}^N e_i}.
\end{align*}
This is tantamount to saying that $(p_n)_{n\ge 1}$ is a Bernoulli sequence with respect to $Q$. Then, \eqref{eq:conv_cantelli} continues to be valid with $Q$ in the place of $P$. Moreover, for each $d$ in $\Omega_1:=\cup_{n\ge 1} S_n$ one has
\[
f_n(d)\to 0 \quad \mbox{ as } \quad n\to\infty.
\]
Thus, since $Q(\Omega_1)=\int_\N Q_n(\cup_{k\ge 1}S_k)\,\gamma(\ddr n)\ge \int_\N Q_n(S_n)\,\gamma(\ddr n)=1$, one sees that $(f_n)_{n\ge 1}$ converges to zero almost surely, and not to $p$. \hfill\qed

\medskip

The phenomenon highlighted in the previous example can be explained in the following terms. There is an instant $n$ beyond which each trial turns into a failure $(=0)$. One is not able to predict when such an instant happens but, according to the definition of $\gamma$, $n$ is viewed as immensely distant. Hence, the probability of events which depend on any finite number of trials are not affected by the instant the sequence becomes a sequences of $0$'s. An analogous example has been provided, more recently, by Ramakrishnan and Sudderth in \cite{rs(88)}. See also the instructive final \textit{Remarks} therein, apropos of the common opinion on finitely additive probabilities.

As pointed out by de Finetti himself, one can change the sequence of $0$'s which, in Example~4.3, follows the last $1$, in such a way that \eqref{eq:conv_cantelli} still continues to be valid and, at the same time, $f_n$ converges or does not: in the former case, it converges to a random variable with a prefixed p.d. function. For the sake of completeness, de~Finetti's paper is here supplemented with a couple of additional examples of that type.
\medskip

\textbf{Example 4.4} \textbf{(a)} Let $\Omega$ be defined as in Example~4.3. Moreover, let $(y_n)_{n\ge 1}$ be the sequence defined by
\[
y_j=\left\{\begin{array}{ll}
0 & \qquad \mbox{ if } \: j\in\{(2n-1)!,(2n-1)!+1,\ldots,(2n)!-1\}\\[7pt]
1 & \qquad \mbox{ if } \: j\in\{(2n)!,(2n)!+1,\ldots,(2n+1)!-1\}
\end{array}
\right.
\]
for $n=1,2,\ldots$. Replace the sets $S_1,S_2,\ldots$ of Example~4.3 by
\[
S_1^*:=\{(1,y_1,y_2,\ldots)\}
\]
and, for $n\ge 2$,
\[
S_n^*:=\{(e_1,\ldots,e_{n-1},1,y_1,y_s,\ldots):\: (e_1,\ldots,e_{n-1})\in\{0,1\}^{n-1}
\}.
\]
Now, for each $n$, set
\[
N_\nu:=n+(2\nu)!-1,\quad M_\nu:=n+(2\nu+1)!-1\qquad\qquad (\nu=1,2,\ldots).
\]
Then, for every sequence in $S_n^*$, one has
\begin{equation}
\left\{\begin{array}{ll}
f_{N_\nu}=\frac{\{(2\nu-1)!+n-1\}\, f_{M_{\nu-1}}}{n+(2\nu)!-1}\,\to\, 0 &\qquad
(\nu\to\infty)\\[8pt]
f_{M_\nu}=\frac{\{(2\nu)!+n-1\}\, f_{N_{\nu}}+2\nu(2\nu)!}{n+(2\nu+1)!-1}\,\to\, 1 &\qquad
(\nu\to\infty).
\end{array}\right.
\label{eq:sequences}
\end{equation}
Define probabilities $Q_1^*,Q_2^*,\ldots$ on $2^\Omega$ according to
\[
Q_1^*(\{(1,y_1,y_2,\ldots\})=1
\]
and, for $n\ge 2$,
\[
Q_n^*(\{e_1,\ldots,e_{n-1},1,y_1,y_2,\ldots\})=p^{\sum_{i=1}^{n-1}e_i}
(1-p)^{n-1-\sum_{i=1}^{n-1}e_i}.
\]
Then, set
\[
Q^*(A):=\int_\N Q_n^*(A)\,\gamma(\ddr n)
\qquad\qquad (A\subset\Omega).
\]
It is easy to prove that $(p_n)_{n\ge 1}$ is a Bernoulli sequence even with respect to $Q^*$. Hence, \eqref{eq:conv_cantelli} is valid even with $Q^*$ in the place of $Q$. On the other hand, in view of \eqref{eq:sequences}, for each sequence in $\Omega_1^*$ one has $\lim\inf f_n=0<\lim\sup f_n=1$, and $Q(\Omega_1^*)=1$.
\smallskip

\noindent \textbf{(b)} Maintain the meaning for $\Omega$, $\gamma$ and $(f_n)_{n\ge 1}$. Denote by $C$ the subset of $\Omega$ on which $(f_n)_{n\ge 1}$ converges and by $H$ a prefixed p.d. function, supported by $[0,1]$. Finally, let $\sigma$ be a $\sigma$--additive probability measure on $\Bcr(\{0,1\}^\infty\cap C)$ such that
\[
\sigma(\{p_1=e_1,\ldots,p_n=e_n\})=\int_{[0,1]}\theta^{\sum_{i=1}^n e_i}
(1-\theta)^{n-\sum_{i=1}^n e_i}\:\ddr H(\theta)
\]
for any $(e_1,\ldots,e_n)\in\{0,1\}^n$ and $n\ge 1$. Now, define
\begin{align*}
  S_n^{**} &:= \{(e_1,\ldots,e_{n-1},1,s):\:(e_1,\ldots,e_{n-1})\in\{0,1\}^{n-1},\: s\in C\} \,\, (n \geq 2) \\
S_1^{**} &:=\{(1,s):\: s\in C\}
\end{align*}
together with the probabilities
\begin{align*}
  Q_1^{**}(\{1\}\times A) &:= \sigma(A)\\
  Q_n^{**}(\{e_1,\ldots,e_{n-1},1\}\times A) &= \beta_{n-1}(\{e_1,\ldots,e_{n-1}\})\,
  \sigma(A)
\end{align*}
for any $A$ in $\Bcr(\{0,1\}^\infty\cap C)$, $(e_1,\ldots,e_{n-1})\in\{0,1\}^{n-1}$, $n\ge 2$ and
\[
\beta_{n-1}(\{e_1,\ldots,e_{n-1}\})=p^{\sum_{i=1}^{n-1}e_i}(1-p)^{n-1-
\sum_{i=1}^{n-1}e_i}.
\]
Finally,
\[
Q^{**}(A):=\int_\N Q_n^{**}(A)\,\gamma(\ddr n)\qquad\qquad (A\subset \Omega).
\]
Once again, $(p_n)_{n\ge 1}$ turns out to be a Bernoulli sequence with respect to $Q^{**}$ and, then, \eqref{eq:conv_cantelli} is valid with $Q^{**}$ in the place of $Q$. On the other hand, in view of de~Finetti's theory of exchangeable sequences, $Q^{**}(\Omega_1^{**})=1$ with $\Omega_1^{**}=\cup_{n\ge 1}S_n^{**}$, and $(f_n)_{n\ge 1}$ converges, on $\Omega_1^{**}$, to a random number whose p.d. function is $H$.
 \hfill \qed
\medskip

The last two groups of examples have an important element in common. They show that there are sequences of random numbers either converging almost surely or oscillating almost surely, with p.d. functions converging weakly in both cases to a p.d. function which differs from the p.d. function of the almost sure limit: The former, in those examples, has a jump invariably equal to 1 at $p$, whereas the latter can be let varying in the class of all p.d. functions. De~Finetti explains why he was interested in investigating into these phenomena in the last section of  \textit{Sui Passaggi}. In the second half of 1929 he was about to deduce the p.d. of the integral of a continuous stochastic f.i.s.i. as limit of the p.d. functions of integral sums which converge pointwise to the integral of interest. At this point, he was assailed by the doubt that such a line of reasoning could be in conflict with his way of thinking of probability, in the sense that the argument could be valid for specific extensions of a prefixed system of finite--dimensional laws, but not in general. The above examples confirmed the reasonableness of his doubt. Here we faithfully follow \cite{df(29c)} and provide the reader with some further insight on this aspect.

Let $t\mapsto X(t)$ be a random function, for $t\ge 0$, with continuous trajectories such that $X(0)\equiv c$. See Section~3. Then, for each $t>0$, one can write
\[
\int_0^t X(u)\,\ddr u=\lim_{n\to\infty}\frac{t}{n}\sum_{h=1}^n
X\left(\frac{th}{n}\right)
\]
where, by the Brunacci--Abel identity,
\[
V_n:=\frac{t}{n}\sum_{h=1}^n
X\left(\frac{th}{n}\right)=ct+\frac{t}{n}\sum_{h=1}^n (n-h+1)\left\{
X\left(\frac{th}{n}\right)-X\left(\frac{t(h-1)}{n}\right)
\right\}.
\]
Then, if $X$ is a f.i.s.i., as proved in \cite{df(29a)}, the characteristic function of the increment $\{X(\frac{th}{n})-X(\frac{t(h-1)}{n})\}$ is given by $\phi(\,\cdot\,)^{t/n}$, with $\phi(\,\cdot\,)$ being the characteristic function of $X(1)$. Then, for the characteristic function $\Psi_{V_n}$ of $V_n$ one gets
\[
\Log\Psi_{V_n}(\xi)=\idr c\xi t+\frac{t}{n}\,\sum_{h=1}^n \Log \phi\left((n-h+1)\frac{t}{n}\,\xi\right) \qquad (\xi\in\R)
\]
where $\Log$ denotes the principal branch of the logarithm. Then, $\lim_{n\to\infty} \Psi_{V_n}$ exists, uniformly on compact intervals, and is given by
\begin{equation}
  \label{eq:repres_id}
  \exp\left\{\idr c\xi t +\frac{1}{\xi}\int_0^{\xi t} \Log(\phi(u))\,\ddr u\right\}.
\end{equation}
After proving that such a limit is a characteristic function, de~Finetti can assert that $V_n$ converges in distribution, but he cannot state that \eqref{eq:repres_id} is the characteristic function of $\lim_n V_n=\int_0^t X(u)\,\ddr u$. In fact, in \cite{brr(07)} it is shown that, in the frame of de~Finetti's theory as summarized in Section~3, $\int_0^1 X(u)\,\ddr u$ can be given any p.d. when, for example, $X$ has the finite--dimensional distributions of the standard Brownian motion.

As a natural development of the previous remarks, de~Finetti tries to find interesting types of convergence of a sequence $(X_n)_{n\ge 1}$ of random numbers to a random number $X$, which entails weak convergence of the p.d. functions of the $X_n$'s to the p.d. of $X$. It is apropos of this question that he introduced the general notion of stochastic uniform convergence (namely convergence in probability) in \cite{df(29c)} and proved that it meets the above property, with respect to weak convergence, in Section~8 of \cite{df(30a)}. See the previous Section~2. Coming back to the problem of the law of the integral of a random function, he concludes \textit{Sui Passaggi} touching on special cases in which one can verify that $(V_n)_{n\ge 1}$ converges in probability to $\int_0^t X(u)\,\ddr u$. This is the case when, for example, $X$ is non--decreasing.

\section{Ensuing correspondence between de Finetti and Fr\'echet}

The correspondence consists of four short open letters published in \textit{Rendiconti del Reale Istituto Lombardo di Scienze e Lettere} and presented in two meetings of the \textit{Istituto}: The first and the second, by Fr\'echet and de~Finetti respectively, were presented on the 3rd of July 1930, whereas the third and the fourth, by Fr\'echet and de~Finetti respectively, on the 20th of November 1930.

In his first note, Fr\'echet begins by saying that \textit{Sui Passaggi} is an interesting paper, which has points of contact with topics he has been dealing with in his recent courses (1929-1930). These were summarized in an article still in press at the beginning of July: that is the paper \cite{f(30c)}, mentioned in Section~2 and published in last issue of \textit{Metron} of that very same year. Fr\'echet agrees with de~Finetti on the fact that countable additivity cannot be deduced from finite additivity, and that the latter constitutes a principle generally accepted as a basis for the theory and the calculus of probabilities. But he has a different opinion about the admissibility of probabilities that are not countably additive. He explains this attitude by referring to the alternative arisen to the founders of the modern theory of measure. They, in spite of the awareness of the impossibility of the problem of measure highlighted by the celebrated Vitali example (see \cite{v(05)}), opted for countable additivity restricted to suitable domains. Fr\'echet holds this way of proceeding up as an example, and proposes introducing the idea of ``\'ev\'enements qui ont une probabilit\'e d\'etermin\'ee et d'autres qui n'en auront pas'', provided that the condition of countable additivity is admitted by definition. Then, he asserts that in such a case, the events considered by de~Finetti in his examples do not have a ``probabilit\'e determin\'ee'', and concludes that it is in this circumstance that the solution to the questions raised by de~Finetti must be sought.

In his answer, de~Finetti respectfully tries to bring the debate down to the real question: To decide if all finitely additive probabilities are admissible, or alternatively if it is necessary to restrict admissible probabilities to the laws that are countably additive. He says he has the sense that several authors, dealing with this subject, consider themselves free to decide according to what suits them best. As an example of this attitude, he mentions the Fr\'echet evocation of an analogy with the theory of measure. Taking his cue from this, de~Finetti says he considers it unjustified to make use of \textit{conventions} to define concepts, like probability, that have a proper meaning, even if possibly open to dispute. Then, the main issue does not consist of making more or less arbitrary conventions on certain properties, but can be traced back to \textit{proving} that certain properties are \textit{necessary}. As mentioned in Section~3, de~Finetti derived the necessity of finite additivity from a coherence principle which, far from being merely conventional, corresponds to a prevailing rational attitude. To say that an event $E$ has probability $p$ either has a more or less common intuitive meaning, or is a perfectly useless sentence. In the first case, if we have a countable family of mutually disjoint events $(E_1,E_2,\ldots)$, with $P(E_n)=0$ for every $n$, are we able to conclude that $E:=\cup_{n\ge 1}E_n$ has probability 0, or, equivalently, that $E$ is invariably practically impossible? It is plain that this is not a convention matter, since the above conclusion has a real conceptual content. It can be false or true, but this must be proved, not assumed by a convention. De~Finetti admits a quite serious difficulty in doing this, due to the lack, until then, of a general theory of probability. He takes the advantage of this circumstance to announce that he has completed a theory, which, starting from a method for assessing probabilities [betting scheme, auth. note], allows one to impeccably deduce the mathematical properties of probability. This is clearly a reference to the theory outlined in the previous Section~3. He recalls that, within such a theory, countable additivity is not a necessary requisite for admissibility and, then, all the examples given in \textit{Sui Passaggi} are perfectly justified and make sense. He concludes with a mention to a couple of consequences, that he considers bizarre, of the adoption of countable additivity as a compulsory principle for admissibility of probability laws. The first consequence relates to the fact that such a principle would forbid one to think of a sequence of mutually disjoint events, forming a partition of the sure event and having probabilities of the same order of magnitude. In other words, that sequence invariably ought to include a finite subset with respect to which the whole of the remaining events would be negligible. The second bizarre consequence of assuming, as a compulsory principle, countable additivity is that one cannot say that the weak limit of a sequence of p.d. functions is always a p.d. function.

The ``official'' reply by Fr\'echet to the previous de~Finetti's arguments was submitted for publication on the 20th of November 1930. The abstract is categoric but the content of the paper is kept to the point more than the first letter. As to the abstract, he writes that de~Finetti's examples are inadmissible, on the basis of the following facts: First, the probabilities of the events considered therein cannot be expressed by real numbers. Second, the probability laws studied therein are inconsistent with empirical experience. Fr\'echet splits his criticism into four points, and de~Finetti answers them in the same order. The final part of this section is accordingly organized into four subsections, each of which summarizes both the Fr\'echet critical remark and the de~Finetti answer pertaining to the point in title of the subsection.
\smallskip

\noindent\textbf{Point 1.} Fr\'echet (F in the sequel) says that de~Finetti (dF in the sequel) has proved, by simple examples, that it is not possible to define a probability for \textit{all} the events while complying, at the same time, with the condition of countable additivity for all classes of mutually disjoint events. Then, he insists on the point, already mentioned in his first letter, that the question can be solved by introducing suitable restrictions on the class of the events equipped with a probability, in such a way that countable additivity is preserved.

dF replies that F has misunderstood his thought. Indeed, he maintains, on the one hand, that it is always possible to comply with the principle of countable additivity and that, on the other hand, all the laws that meet the principle of coherence are admissible. The first part of this argument is not in contrast with the Vitali theorem, which simply excludes that a countably additive law might give equal probabilities to all superposable subsets of $[0,1]$. Then, to overcome the drawback, it is not necessary to restrict the class of the admissible events. Moreover, such a strategy is not so much as sufficient. To see this, reconsider Example~4.1 with $\Omega=(0,1]$ and $x_1=1$, $x_2=1/2$, $x_3=1/3$,~... and form the partition of $(0,1]$ in the intervals given by $(1/2,1]$, $(1/3,1/2]$,~..., $(1/(n+1),1/n]$,~...~. Then, since $P(\{x_1,\ldots,x_n,\ldots\})=1$, one gets $P((0,1])=1$, but $P((1/(n+1),1/n])=0$ for every $n$. In other words, the ``drawback'' of countable partitions of $\Omega$, into elements having zero probabilities, can occur even if these elements are intervals, which, according to F, constitute  the most typical example of ``probabilizable events''. It is at this point that dF drops the artificial examples suggested by the theory of measure, to explain how more concrete situations, pertaining to the study of functions with random increments, are open to the same objections as those raised in the previous more schematic examples. To this end, dF mentions two cases. One of them involves the notion of \textit{conglomerability} and, hence, goes beyond the scope of the present paper. The original paper in which that phenomenon had been noticed is \cite{df(30g)}. The other case is concerned with the already mentioned de~Finetti's theorem on the nowhere differentiability of the trajectories of the Brownian motion. See Section~2. He expresses his regret that, in view of his criticism on the role of countable additivity, the aforementioned theorem must be reformulated in a weaker form, i.e.: Let $\ep$ and $M$ be strictly positive numbers. Then, the probability that $[0,1]$ includes any interval of length greater then $\ep$, for which one gets $|X(t_2)-X(t_1)|<M(t_2-t_1)$, is zero. Since the usual formulation, obtained as $M\to\infty$ and $\ep\to 0$, would be very important, if valid in general, dF admits he would be very happy of the existence of any reasonable argument that persuades him to share the common idea that probability laws are continuous. He concludes expressing his skepticism towards the solution devised by F: ``Even supposing that there are events for which the doubts about countable additivity turn out to be groundless, how could I recognize them in practical situations of the same type as that just now described.''
\smallskip

\noindent\textbf{Point 2.} This corresponds to the first point raised in the above-mentioned abstract. F seems to admit that restricting the class of the ``probabilizable events'' does not serve the purpose to explain the antinomy stressed by de~Finetti's examples. In order to solve the issue, he proposes to consider probabilities expressed in terms of actual infinitesimals, say $\bm{\ep}$. He supposes that, in such a way, one succeeds in writing $\bm{\ep}\cdot\bm{\omega}=1$. dF observes that the adoption of ``new numbers'' of the type of $\bm{\ep}$  is not inconsistent with his own way of thinking. In fact, he recalls he has already made use of those numbers in \cite{df(28)}. But, unlike F, he gets to the conclusion there is no contradiction between admissibility of infinitesimals and finite additivity. Indeed, the probability of the union of any finite number of infinitesimal events is infinitesimal and no limit process could lead to conclude that the probability of the union of all the events is $1$.
\smallskip

\noindent\textbf{Point 3.} It corresponds to the second point briefly described in the abstract. F maintains that probabilities like those of Example~4.1 do not appear when probabilities are based on frequencies. Moreover, he is skeptical about the success of dF's scientific plan, since none of the definitions of probability proposed until then had met with general approval. So, he refers the reader to the final part of the L\'evy monograph \cite{l(25)}, where countable additivity is, according to his opinion, justified. Apropos of the first assertion, he claims he is able to prove it in the following terms. Consider a random phenomenon with a countable set of elementary possible outcomes, say $a_1,a_2,\ldots$, and assume one can conceive an indefinitely extendable sequence of trials of that phenomenon. Let $f_k^{(n)}=r_k/n$ be the frequency of $a_k$ in the first $n$ trials. Then, all but a finite number of the $f_k^{(n)}$'s are zero and $f_1^{(n)}+f_2^{(n)}+\,\cdots=1.$ It is palpably clear that such an equality will be valid for any $n$. According to the empirical interpretation, as explained in \cite{fh(24)}, $f_k^{(n)}$ represents an ``experimental measure'' of a probability $p_k$, for every $k$, as $n$ increases. F claims that these remarks are sufficient to conclude that $p_1+p_2+\,\cdots\,=1$.

In his reply, dF stresses once again that properties like countable additivity be proved and not established in the form of conventions. Even if he admits the existence of the difficulties mentioned by F, he considers them more extrinsic rather than intrinsic. He says that, to clear this hurdle, it would suffice that every author would give his own proof, based on his own definition of probability. He acknowledges that this is what F has tried to do, starting from his empirical interpretation of probability. But, despite the inadequacy of the F proof, dF shows that the F argument, once made more precise, can become an excellent point in favour of the dF thesis. As to the inadequacy, assuming that each $p_k$ is the (usual) limit, as $n\to\infty$, of $(f_k^{(n)})_{n\ge 1}$, one can write $1=\lim_{n\to\infty}\sum_{k} f_k^{(n)}$ but, in general, the exchange of $\lim$ with $\sum$ is not valid. Hence, one can just say, in general, that $1\ge \sum_k p_k$. In the footnote (2) of page~256, dF notes that F had pointed out, in personal correspondence, that he didn't mean to speak of $p_k$ as a limit in a mathematical sense. Once taken note of this detail, dF resorts to a different argument free from the criticism of being just in an abstract mathematical form.  He sets himself the objective of studying the expected behaviour of the $f_k^{(n)}$'s, as $n$ goes to infinity, to show that \textit{there are probability laws with respect to which it would be illusory to expect that the $f_k^{(n)}$'s converge to numbers $p_k$ such that $\sum p_k=1$.} Therefore, in the same way, it would be illusory to hope to prove that the property of countable additivity may be derived, in general, from the analogous property valid for frequencies, understood as \textit{empirical} estimates of probabilities. Here is the example proposed by dF along with a few further details.
\smallskip

\textbf{Example 5.1.} Let $S$ be a countably infinite set, say $S:=\{a_1,\ldots,a_n,\ldots\}$, and $\Omega=S^\infty$. Define $\xi_1,\xi_2,\ldots$ to be the coordinate random variables of $\Omega$. Now, for any strictly positive integer $N$, set $S_N:=\{a_1,\ldots,a_N\}$ and, for every $A\subset S^n$, define
\[
P_N^{(n)}(A)=\frac{\sharp(A\cap S_N^n)}{N^n}.
\]
It is plain that $(S^n,2^{S^n},P_N^{(n)})_{n\ge 1}$ is a consistent system of probabilities. Then, there exists a probability $P_N$ on $(\Omega, 2^\Omega)$ such that $P_N(A\times S^\infty)=P_N^{(n)}(A)$ holds true for every $A\subset S^n$, $n=1,2,\ldots$. Finally, with the same $\gamma$ as in Example~4.3, put
\[
P(C)=\int_\N P_N(C)\,\gamma(\ddr N) \qquad\qquad (C\subset\Omega).
\]
Hence, $P$ is a probability on $2^\Omega$ such that $P\{\xi_n=a_k\}=0$ for every $k$ and $n$. Indeed, for every $N\ge k$, one has
\[
P_N\{\xi_n=a_k\}
=P_N^{(n)}(S^{n-1}\times\{a_k\})=\frac{\sharp((S^{n-1}\times\{a_k\})\cap S_N^n)}{N^n}
=\frac{N^{n-1}}{N^n}
\]
and, then, $P\{\xi_n=a_k\}=\int_{\{N\ge k\}}N^{-1}\,\gamma(\ddr N)=0$. Moreover, for any $n\ge 2$,
\begin{equation}
  \label{eq:different_coord}
  P\{\xi_1\ne \xi_2\ne\,\cdots\,\ne \xi_n\}=1
\end{equation}
Indeed, if $N>n$, one gets
\begin{align*}
  P_N(\{\xi_1\ne \xi_2\ne\,\cdots\,\ne \xi_n\}\times S^\infty)
  &=P_N^{(n)}\{\xi_1\ne \xi_2\ne\,\cdots\,\ne \xi_n\}\\[7pt]
  &=\frac{N(N-1)\,\cdots\,(N-n+1)}{N^n}
\end{align*}
and, since
\[
P\{\xi_1\ne \xi_2\ne\,\cdots\,\ne \xi_n\}=\int_{\{N\ge n\}}\frac{N(N-1)\,\cdots\,(N-n+1)}{N^n} \gamma(\ddr N)
\]
the equality in \eqref{eq:different_coord} follows.
Now, if $(p_k)_{k\ge 1}$ is any sequence such that $p_k\ge 0$ for every $k$ and $\sum_{k} p_k=1$, there is $\bar k$ such that $p_{\bar k}=\max\{p_1,p_2,\ldots\}>0$. Then, there is a contradiction between the adoption of $P$ as a probability on $2^\Omega$ and the assumption that $f_k^{(n)}:=\sum_{j=1}^n\indic_{\{a_k\}}(\xi_j)/n$ converges, in some sense, to $p_k$ for every $k$, as $n\to\infty$. Indeed, with respect to $P$, for each $n$ it is practically sure~-- in view of \eqref{eq:different_coord}~-- that there are $n$ indices, say $k_1,\ldots,k_n$, for which $f_{k_i}^{(n)}=1/n$ $(i=1,\ldots,n).$ So, one gets
\begin{equation}
  \label{eq:freq}
  P\left\{f_{\bar k}^{(n)}\le\frac{p_{\bar k}}{M}\right\}=1
\end{equation}
for any $M,n\in\N$ and $n>M/p_{\bar k}$, thus contradicting any reasonable definition of convergence of $f_{\bar k}^{(n)}$ to $p_{\bar k}$, as $n\to\infty$. \hfill\qed
\smallskip

The above example shows that the condition of countable additivity does not follow from the fact that $\sum_k f_k^{(n)}=1$ holds true in any case, combined with the assumption that frequencies approach probabilities as $n$ increases to infinity. Indeed, \eqref{eq:freq} and the arbitrariness of $(p_k)_{k\ge 1}$ show, once again, the inadequacy of the use of frequencies to prove countable additivity. The sense of this statement can be strengthened even further when $P_N$ is defined to be the Kolmogorov extension of the $P_N^{(n)}$'s to the smallest $\sigma$--alegbra containing all the sets $\xi_m^{-1}(A)$ for all $m$ and all $A\subset S$. It is easy to check that the $\xi_n$'s turn out to be independent and identically distributed, with uniform distribution on $S_N$, with respect to any extension of the $P_{N}^{(n)}$'s. Moreover, if $P_N$ corresponds to the Kolmogorov extension, then the strong law of large numbers yields
\[
P_N\left\{\lim_n f_k^{(n)}=\frac{1}{N} \mbox{ for } k=1,\ldots,N,
\quad \lim_n f_k^{(n)}=0 \mbox{ for } k\ge N+1\right\}=1
\]
and, with the same $(p_n)_{n\ge 1}$ as in Example~5.1 and for every $M$ in $\N$ one obtains
\begin{align*}
  P\left\{\sum_k \overline{\lim}_{n\to\infty} f_k^{(n)}\le \frac{1}{M}\right\}
  & \ge P\left\{\overline{\lim}_{n\to\infty} f_k^{(n)}\le \frac{p_k}{M}
  \mbox{ for every } k\right\}\\
  &=\int_\N P_N\left\{\overline{\lim}_{n\to\infty} f_k^{(n)}\le \frac{p_k}{M}
  \mbox{ for every } k\right\}\:\gamma(\ddr N)\\
  &=1.
\end{align*}
Then, with respect to this particular $P$, we are practically sure that frequencies converge and that the sum of the series is in $[0,\ep)$ for every $\ep>0$. \hfill \qed
\smallskip

\noindent\textbf{Point 4.} The last objection raised by F is about a seeming slip made by dF, in the previous part of the correspondence, apropos of the nature of the weak limit of a sequence of p.d. functions. As already recalled in this very same section, dF found it bizarre that, in the common approach based on countable additivity, such limit was not necessarily a p.d. function. F maintains that this phenomenon is not warded off by the adoption of dF point of view. In support of this statement, F gives the example of the p.d. functions $F_n(x)=(\indic_{(-n,n]}/2 +\indic_{(n,\infty)})(x)$, $x\in\R$ and $n=1,2,\ldots$. It is easy for dF to prove that the F argument is ineffective since the function $F\equiv 1/2$, limit of $(F_n)_{n\ge 1}$, can be viewed, within the finitely additive frame, as a p.d. function.

\section{Final remarks} The paper we have annotated in the previous sections cannot be counted  among de~Finetti's most important works. Nonetheless, it represents a direct evidence of an extremely interesting stage of his scientific career. It is the stage of the mathematical formulation of his subjectivistic conception of probability, and of the consequent conclusion that the only general restriction on the class of the admissible probability laws is given by the coherence principle. Hence, probabilities on algebras of events must be additive, but not necessarily countably additive. This, on the one hand, led him to revise the value, in terms of their generality, of a few of his own previous theorems that were proved, in part, under the assumption of continuity of probability laws. His critique, on the other hand, did not even spare one of the most renowned achievements of the theory of probability, i.e. the strong law of large numbers. By resorting to enlightening examples, in \textit{Sui Passaggi} he succeeds in enhancing some crucial differences between the two viewpoints taken into consideration therein. Besides, the discussion with Fr\'echet gives de~Finetti a chance to provide fresh and deep explanations about his stance. They continue to be of great interest and useful since the arguments, still put up against the adoption of de~Finetti's theory, do not basically differ from those used by Fr\'echet.

\end{document}